\documentclass{amsart}

\RequirePackage{amsmath}
\RequirePackage{amssymb}
\RequirePackage{xy}

\xyoption{all}
\CompileMatrices

\theoremstyle{plain}
\newtheorem{theorem}{Theorem}[section]
\newtheorem{proposition}[theorem]{Proposition}
\newtheorem{lemma}[theorem]{Lemma}

\newtheorem{conjecture}[theorem]{Conjecture}
\newtheorem*{mthm}{Theorem}
\newtheorem*{mprop}{Proposition}
\newtheorem*{mconj}{Conjecture}

\theoremstyle{definition}

\newtheorem{remark}[theorem]{Remark}

\numberwithin{equation}{section}

\newcommand{\A}{\mathbf{A}}
\renewcommand{\a}{\mathfrak{a}}

\newcommand{\Fp}{{\mathbf F}_{p}}
\newcommand{\Fpbar}{\bar{\mathbf F}_{p}}

\renewcommand{\O}{\mathcal{O}}
\newcommand{\p}{\mathfrak{p}}
\renewcommand{\P}{\mathcal{P}}
\newcommand{\PP}{\mathfrak{P}}
\newcommand{\Q}{\mathbf{Q}}
\newcommand{\Z}{\mathbf{Z}}

\newcommand{\Nm}{\mathbf{N}}

\newcommand{\pbar}{\bar{\p}}
\newcommand{\pib}{\bar{\pi}}
\newcommand{\Qbar}{\bar{\Q}}
\newcommand{\Ql}{\Q_{\ell}}
\newcommand{\Qp}{\Q_{p}}
\newcommand{\sd}{\sqrt{-d}}
\newcommand{\vep}{\varepsilon}

\newcommand{\ls}[2]{\left(\frac{#1}{#2}\right)}

\newcommand{\ts}{\textstyle}

\newcommand{\amod}[1]{\,(\bmod{\,#1})}

\newcommand{\oh}{\ts \frac{1}{2}}
\renewcommand{\th}{\text{th}}
\newcommand{\surj}{\twoheadrightarrow}

\DeclareMathOperator{\Frob}{Frob}
\DeclareMathOperator{\Gal}{Gal}
\DeclareMathOperator{\GL}{GL}
\DeclareMathOperator{\tr}{tr}
\DeclareMathOperator{\cm}{cm}

\begin{document}

\title{Power residues of Fourier coefficients of modular forms}
\author{Tom Weston}

\maketitle

Let $f = \sum a_{n}q^{n}$ be a newform with rational Fourier coefficients.
For an integer $m \geq 2$
let $\delta_{m}(f)$ denote the relative density of the set
\begin{equation} \label{eq:set}
\bigl\{ p \equiv 1 \amod{m} \,;\, a_{p} \text{~is a non-zero~} m^{\th} 
\text{~power modulo~} p \bigr\}
\end{equation}
inside the set of primes $p \equiv 1 \amod{m}$
such that $a_{p}$ is non-zero modulo $p$.
Based on computations with various newforms we make the following
conjecture.

\begin{mconj}
If $f$ does not have complex multiplication, then
$\delta_{m}(f) = \frac{1}{m}$.
\end{mconj}

In fact, we suspect that much more is true: we conjecture that
this relative density does not change after restriction to any set
of primes defined by a Cebatorev-style Frobenius condition;
that is, we expect that the sets (\ref{eq:set}) 
yield sets of primes of positive density which are quite different from
those sets determined by Galois theoretic conditions.  See
Conjecture~\ref{conj2} for a precise statement.

For CM-forms $f$ and certain values of $m$ the set (\ref{eq:set})
is in fact defined by Galois theoretic conditions.
We use this to prove the following result;
see Theorems~\ref{thm:sqrs},~\ref{thm:cube}, and~\ref{prop:high}
for precise statements.  (In particular,
Theorem~\ref{thm:sqrs} includes the cases $d=1,3$ when $m=2$.)

\begin{mthm}
Let $K=\Q(\sqrt{-d})$ be an imaginary quadratic field with odd class
number $h$.  Fix $k \geq 2$, $k \equiv 1 \amod{h}$, 
and let $S_{k}^{K-\cm}(\Q)$ denote the
set of newforms of weight $k$ with rational coefficients and
complex multiplication by $K$.
\begin{enumerate}
\item Assume that $d \neq 1,3$ and $k$ even.  Then $\delta_{2}(f)=\frac{1}{2}$
for all but two $f \in S_{k}^{K-\cm}(\Q)$; for these two forms,
$\delta_{2}(f)$ equals $\frac{1}{4}$ or $\frac{3}{4}$.
\item Assume that $d \neq 1,3$ and $k$ odd.  Then $\delta_{2}(f)=\frac{3}{4}$
for all but two $f \in S_{k}^{K-\cm}(\Q)$; for these two forms,
$\delta_{2}(f)$ equals $\frac{1}{2}$ or $1$.
\item Assume that $d=3$.  Then $\delta_{3}(f)$ equals either 
$\frac{5}{9}$ or $1$ for $f \in S_{k}^{K-\cm}(\Q)$.
\item Assume that $d > 3$.  Then for any $m | k-1$ we have
$\delta_{m}(f) = \frac{3}{4}$ for all but finitely many
forms $f \in S_{k}^{K-\cm}(\Q)$; for these exceptional forms,
$\delta_{m}(f)$ equals either $\frac{1}{2}$ or $1$.
\end{enumerate}
\end{mthm}

We fully expect that similar results can be proven for imaginary
quadratic fields with even class number, but we do not consider that
case here.

The original motivation for this work was the following question of 
Ramakrishna: if $E$ is an elliptic curve over $\Q$, are the Fourier
coefficients $a_{p}(E) := p + 1 - \#E(\Fp)$ cubes for infinitely many
$p$?  This question in turn was motivated by the following observation
of Ramakrishna (in the case $m=3$), which we discuss in detail in Section 2.3.

\begin{mprop} 
Let $K$ be an imaginary quadratic field of class number one and let
$m$ be a prime relatively prime to $\# \O_{K}^{\times}$.  Let
$p \equiv 1 \amod{m}$ be a prime greater than $3$ 
which splits in $K/\Q$ and let
$K_{p}^{m}$ be the maximal abelian $m$-extension of $K$ which is
unramified away from $p$.  Then $p$ has inertial degree one in
$K_{p}^{m}/\Q$ if and only if $a_{p}(E)$ is an $m^{\th}$ power modulo
$p$, where
$E$ is any rational elliptic curve with complex multiplication by $K$
and good reduction at $p$.
\end{mprop}

In fact, one can state the above criterion in terms of modular forms
of higher weight as well; see Remark~\ref{rmk:hwt}.  
Unfortunately, our methods do not yield a single case in which
we can show that the criterion of the proposition is satisfied for
infinitely many $p$.

In Section 1 we formulate our precise conjectures and discuss the
numerical evidence.  In Section 2 we begin the study of the CM-case,
relating the power residues of Fourier coefficients of
CM-forms to power residues of the associated Hecke character.
In Section 3 we combine this formula with the quadratic reciprocity law
to compute $\delta_{2}(f)$ for CM-forms $f$.  We consider the cases $d=m=3$
and $m|k-1$ for $f$ of weight $k$ in Section 4.

We remark that the methods of this paper are essentially elementary.
We suspect that it will require much deeper methods 
to make any progress on our conjectures in the non-CM case.

It is a pleasure to thank Ravi Ramakrishna for 
suggesting this problem and for all of the encouragement and
insight he provided.  I would also like to thank 
Rob Benedetto, Ken Ribet, and Siman Wong for
helpful conversations.

\subsection*{Notation}

For $m \geq 2$ and $\p$ a prime of a number field $K$ with
residue field of order $q$ congruent to $1$ modulo $m$ we write
$\ls{\cdot}{\p}_{m}$
for the $m^{\th}$ power residue symbol modulo $\p$; thus for
$\alpha \in \O_{K}$, $\alpha \notin \p$, 
we have $\ls{\alpha}{\p}_{m} \in \mu_{m}$,
$$\ls{\alpha}{\p}_{m} \equiv \alpha^{(q-1)/m} \pmod{\p},$$
and $\ls{\alpha}{\p}_{m} = 1$ if and only if $\alpha$ is a non-zero
$m^{\th}$ power residue modulo $\p$.
We simply write $\ls{\alpha}{\p}$ in the case $m=2$.

By a {\it Hecke character} over a number field $K$ we mean a continuous
homomorphism
$$\chi : \A_{K}^{\times}/K^{\times} \to \Qbar^{\times}$$
with $\A_{K}$ the adeles of $K$.  We say that $\chi$ is {\it unramified}
at a prime $\p$ of $\O_{K}$ if $\chi(\O_{K_{p}}^{\times})=1$ (where
we embed the completion $K_{\p}^{\times}$ of $K$ at $\p$ 
into $\A_{K}^{\times}$ in the obvious manner).
In this case we write $\chi(\p)$ for the value of $\chi$ on any
uniformizer of $K_{\p}$ and we extend $\chi$ to a character on all unramified
fractional ideals in the obvious way.

If $\P$ is a set of primes, by the {\it density} of $\P$ we always mean
the Dirichlet density, although all results in this paper remain true for
natural density as well.  For $\alpha \in \Q^{\times}$ we write
$D(\alpha)$ for the discriminant of $\Q(\sqrt{\alpha})$ over $\Q$.

\section{Conjectures}

\subsection{Galois representations}

Consider a Galois representation
$$\rho : \Gal(\Qbar/\Q) \to \GL_{n}\Ql.$$
Suppose that $\rho$ is {\it motivic} in the sense that there is a
smooth, projective variety $X$ over $\Q$ and a projector $p$
in the ring of algebraic correspondences on $X$ such that $\rho$ is
the Galois representation on
$$p_{*} H_{\text{\'et}}^{i}\bigl(X_{\Qbar},\Ql(j)\bigr)$$
for some $i,j$.  (According to the Fontaine--Mazur conjectures
\cite{FM} it should be equivalent to suppose that $\rho$ is finitely
ramified and potentially semistable at $\ell$.)  It follows from
Deligne's proof of the Weil conjectures \cite{Deligne} that 
$$a_{p}(\rho) := \tr \rho(\Frob_{p})$$
is a rational integer for almost all primes $p$, independent of the
choice of arithmetic Frobenius element $\Frob_{p}$.

For an integer $m \geq 2$, define
$$\delta_{m}(\rho) = \frac{\text{density of~}
\left\{ p \equiv 1 \amod{m} \,;\, \ls{a_{p}(\rho)}{p}_{m} = 1 \right\}}
{\text{density of~} \{p \equiv 1 \amod{m} \,;\, a_{p}(\rho)
\not\equiv 0 \amod{p} \}}$$
if these densities exist.
With Ramakrishna, we make the following conjecture.

\begin{conjecture} \label{conj1}
Let $\rho$ be a motivic Galois representation as above.
If the image of $\rho$ is open, then
$$\delta_{m}(\rho) = \frac{1}{m}.$$
\end{conjecture}

We remark that Conjecture~\ref{conj1} is certainly false for many
infinitely ramified Galois representations, as by \cite{Ravi} it is
possible to control the Frobenius polynomials (and thus the
$a_{p}(\rho)$) at a set of primes of density one.

Note that when $\rho$ has open image one expects the
set of $p$ dividing $a_{p}(\rho)$ to have density
zero, so that $\delta_{m}(\rho)$ should be an absolute density in this case.

Of course, one could attempt to formulate an analogous conjecture for
Galois representations of number fields taking values in larger $\ell$-adic
fields.  We have not attempted to do so here; we only point out that
the density $\delta_{m}(\rho)$ is certainly dependent on the choice of
coefficient field, so that the conjectures must be somewhat more complicated
in the general case.

\subsection{Fourier coefficients of modular forms}

The motivation for Conjecture~\ref{conj1} comes from numerical
investigations with Fourier coefficients of modular forms.
Let $f = \sum a_{n}q^{n}$ be a newform of weight $k \geq 2$ and level
$N$ with rational Fourier coefficients.  For any prime $\ell$
Deligne has constructed a motivic Galois representation
$$\rho_{f} : \Gal(\Qbar/\Q) \to \GL_{2} \Ql$$
with the property that $\tr \rho_{f}(\Frob_{p}) = a_{p}$ for 
$p \nmid N\ell$; see \cite{Scholl}, for example.
By \cite[Theorem 5.7]{Ribet}
the image of $\rho_{f}$ is open provided that $f$ is not of
CM-type.

We briefly discuss the numerical evidence for Conjecture~\ref{conj1}
for the representations $\rho_{f}$; we simply write $\delta_{m}(f)$ for
$\delta_{m}(\rho_{f})$ in this case, or $\delta_{m}(E)$ if $f$ is
the newform corresponding to a rational elliptic curve $E$.
For integers $P_{1} < P_{2}$ define
$$\delta_{m}(f;P_{1},P_{2}) = \frac{\#
\left \{ p \equiv 1 \amod{m} \,;\, \ls{a_{p}}{p}_{m} = 1,\, P_{1} \leq 
p \leq P_{2}\right\}}
{\# \{p \equiv 1 \amod{m} \,;\, a_{p} \not\equiv 0 \amod{p},\, P_{1} \leq p 
\leq P_{2}\}}.$$
We have computed:
\begin{enumerate}
\item $\delta_{m}(E;10^{8},2\!\cdot\! 10^8)$ 
for $m \leq 10$ and various rational elliptic curves of
small conductor.
\item $\delta_{m}(f;1,1000)$ for $m \leq 10$ and various modular forms
$f$ of weight at least $3$ contained in the tables \cite{Stein}.
\item $\delta_{m}(\Delta;10^{6},2\!\cdot\! 10^{6})$ for $m \leq 10$ and the
modular form $\Delta = \sum \tau(n)q^{n}$ with $\tau(n)$ the Ramanujan
$\tau$-function.
\end{enumerate}
In each case we obtained results consistent with 
Conjecture~\ref{conj1}.  To give a single example, we report the information
obtained for the elliptic curve $E=X_{0}(11)$.
The data reported below is for the set of primes
$$\{p \equiv 1 \amod{m} \,;\, 10^8 \leq p \leq 2 \!\cdot 10^8\! \}.$$

\vspace{0.1cm}

\begin{center}
\begin{tabular}{c|c|c|c}
$n$ &  $\# p \,\,\text{s.t.\ } \!\ls{a_{p}}{p}_{m}=1$ & 
$\#p \,\,\text{s.t.\ } \!p \nmid a_{p}$ & 
$\delta_{m}(E;10^8,2 \!\cdot \!10^8)$ 
\\ 
\hline
2 & 2662953 & 5317482 & 0.5008 \\
3 & 888792 & 2658461 & 0.3343 \\
4 & 667722 & 2658316 & 0.2512 \\
5 & 266666 & 1329469 & 0.2006 \\
6 & 446913 & 2658461 & 0.1681 \\
7 & 127203 & 886591 & 0.1435 \\
8 & 168427 & 1329053 & 0.1267 \\
9 & 99178 & 886298 & 0.1119 \\
10 & 133116 & 1329469 & 0.1001 
\end{tabular}
\end{center}

\vspace{0.2cm}

S.\ Wong has pointed out that 
in most cases the data appears to be slightly
biased, so that the approximate densities are usually larger than
$\frac{1}{m}$.  For the case of elliptic curves this may be at least
partially explained by the fact that $|a_{p}| \leq 2\sqrt{p}$ and
such small values are slightly more likely to be power residues.
In any event, Conjecture~\ref{conj1} asserts that in the
limit these biases disappear.

\subsection{Cebatorev sets}

We have also computed densities as above for sets of primes satisfying
additional congruence conditions (that is, with prescribed splitting in
cyclotomic fields) and with specified inertial degrees in 
splitting fields of various cubic polynomials.
These computations suggest a stronger statement than
Conjecture~\ref{conj1}.
For $\rho$ as above and for $\P$
a set of rational primes of positive density, we define
$$\delta_{m}(\rho;\P) = 
\frac{\text{density of~}
\left\{ p \equiv 1 \amod{m} \,;\, \ls{a_{p}(\rho)}{p}_{m} = 1 \right\} \cap \P}
{\text{density of~} \{p \equiv 1 \amod{m} \,;\, a_{p} \not\equiv 0 \amod{p}\} 
\cap \P}$$
if these densities exist.  We say that such a set
$\P$ is a {\it Cebatorev set} if there is a finite Galois
extension $K$ of $\Q$ and a subset $S \subseteq \Gal(K/\Q)$, stable
under conjugation, such that, up to finite sets, $\P$ is the
set of primes $p$ with $\Gal(K/\Q)$-Frobenius lying in $S$.

\begin{conjecture} \label{conj2}
Let $\rho$ be a motivic Galois representation with open image as above.  Then
for any Cebatorev set $\P$,
$$\delta_{m}(\rho;\P) = \frac{1}{m}.$$
\end{conjecture}

Conjecture~\ref{conj2} essentially asserts that the sets
\begin{equation} \label{eq:pset}
\left\{p \equiv 1 \amod{m} \,;\, \ls{a_{p}(\rho)}{p}_{m} = 1 \right\}
\end{equation}
can not be described in terms of Cebatorev sets.  If
Conjecture~\ref{conj2} holds, it would thus yield an entirely new
collection of naturally occurring sets of primes of positive density.

\subsection{CM-representations}

The case where $\rho$ has smaller image in $\GL_{n}\Ql$ does not
appear to admit a uniform statement such as Conjecture~\ref{conj1}.  
We do not speculate on the
form of a general conjecture.
However, in certain cases we expect that the analogous result holds.

\begin{conjecture} \label{conj3}
Let $f$ be a rational newform of weight $k \geq 2$ with complex multiplication
by the field $\Q(\sqrt{-d})$.
Assume that $d > 3$.
If $m$ is relatively prime to $2(k-1)$, then
$$\delta_{m}(f) = \frac{1}{m}.$$
\end{conjecture}

In the remainder of the paper we consider the exceptional cases $m=2$ and
$m|k-1$ (and $m=d=3$).  In these cases we will see that, in contrast
to Conjecture~\ref{conj2},
(\ref{eq:pset}) is a Cebatorev set and thus that we can compute its density.

\section{Modular forms of CM-type with rational Fourier coefficients}

For the remainder of the paper we fix an imaginary quadratic field
$K$ with ring of integers $\O$.
We write $d$ for the unique squarefree integer such that
$K=\Q(\sd)$, $D$ for the discriminant of $K$, and
$w$ for the order of $\O^{\times}$.
For $\alpha \in K$, we write $\bar{\alpha}$ for the conjugate of $\alpha$.

Let $H$ denote the Hilbert class field of $K$, so that $[H:K]$ equals
the class number $h$ of $K$.
We always assume that $h$ is odd, so that
either $d \in \{1,2\}$ or $d \equiv 3 \amod{4}$.
For later use we set $\vep = 2$ (resp.\ $\vep = -1$, resp.\ $\vep = 1$)
for $d=1$ (resp.\ $d \equiv 3 \amod{8}$, resp.\ $d \equiv 7 \amod{8}$ or
$d=2$).

\subsection{Hecke characters}

By
\cite[Section 11.2]{Gross}, for $d > 3$ there is a unique Hecke character
$$\psi' : \A_{H}^{\times}/H^{\times} \to K^{\times},$$
unramified away from $D$, with the property that
for any prime $\PP$ of $\O_{H}$ (relatively prime to $D$)
$\psi'(\PP)$ is
the unique generator of the principal ideal
$N_{H/K}\PP$ which is a square modulo $\sqrt{-d}$.
(Note that $w=2$ and $d \equiv 3 \amod{4}$, so that there is indeed a unique
such generator.)
Since $\psi' \circ \sigma$ has the same property for any
$\sigma \in \Gal(H/K)$, the Hecke character $\psi'$ is invariant under
the action of $\Gal(H/K)$ on $\A_{H}^{\times}/H^{\times}$.

Define a Hecke character $\psi : \A_{K}/K^{\times} \to K^{\times}$
as the composition
$$\A_{K}^{\times}/K^{\times} \to \A_{H}^{\times}/H^{\times}
\overset{\psi'}{\longrightarrow} K^{\times};$$
it follows from Lemma~\ref{lemma:gen} below that $\psi$ has
infinity type $(h,0)$ (in the sense of \cite[Section II.1.1]{deShalit}).

\begin{lemma} \label{lemma:gen}
Let $\p$ be a prime ideal of $\O$ relatively prime to $D$.
Then
$\psi(\p)$ is the unique generator of $\p^{h}$ which is
a square modulo $\sqrt{-d}$.
\end{lemma}
\begin{proof}
Let $\varpi_{\p} \in \A_{K}^{\times}$ be an idele which is trivial away
from $\p$ and which is a uniformizer at $\p$.  Then the image of
$\varpi_{\p}$ in $\A_{H}^{\times}$ can be written as
$\varpi_{\PP_{1}} \cdots \varpi_{\PP_{g}}$ for analogously defined
ideles $\varpi_{\PP_{i}}$; here $\p\O_{H} = \PP_{1} \cdots \PP_{g}$ is the
decomposition of $\p$ into primes of $\O_{H}$.
Since $\psi'$ is $\Gal(H/K)$-invariant it follows that
$\psi(\p) = \psi'(\PP_{1}) \cdots
\psi'(\PP_{g}) = \psi'(\PP_{1})^{g}$.
As $\psi'(\PP_{1})\O_{K} = N_{H/K}\PP_{1} = \p^{h/g}$, we conclude that
$\psi(\p)\O_{K} = \p^{h}$.  Each $\psi'(\PP_{i})$ is a square modulo
$\sqrt{-d}$, so that the same is true of $\psi(\p)$; thus
$\psi(\p)$ is the unique generator of $\p^{h}$ which is a square
modulo $\sqrt{-d}$, as claimed.
\end{proof}

We define the Hecke character $\psi : \A_{K}^{\times}/K^{\times} \to
K^{\times}$ for $d=1$ (resp.\ $d=2$, resp.\
$d=3$) as
the Hecke character over $K$ associated (as in \cite[Theorem 9.2]{Silverman2})
to the $\Q$-isogeny class of
elliptic curve $32A$ (resp.\ $256D$, resp.\ $27A$) of \cite{Cremona}.

\begin{lemma} \label{lemma:gen2}
Let $\p$ be a prime ideal of $\O$ relatively prime to $D$.
Then
for $d=1$ (resp.\ $d=3$)
$\psi(\p)$ is the unique generator of $\p$ which is
congruent to $1$ modulo $2+2i$ (resp.\ modulo $3$).
For $d=2$, $\psi(\p)$ is the unique generator of $\p$ which
is congruent to one of
\begin{equation} \label{eq:2set}
\bigl\{1,3,5+\sqrt{-2},7+\sqrt{-2},5+2\sqrt{-2},7+2\sqrt{-2},5+3\sqrt{-2},
7+3\sqrt{-2}\bigr\}
\end{equation}
modulo $4\sqrt{-2}$.
\end{lemma}
\begin{proof}
This is well-known for $d=1,3$; see for example
\cite[Example II.10.6 and Exercise II..34]{Silverman2}.
For $d=2$, by \cite[Proposition 10.4]{Silverman2} $\psi(\p)$ is a generator
of $\p$ for all $\p \neq \sqrt{-2}\O$.  Since $\psi$ has conductor
dividing $\sqrt{256}=16$, to determine $\psi$ it suffices to
determine the sign of this generator for the principal ideals
generated by representatives for all classes in $(\O/16\O)^{\times}$.
This is straightforward via \cite[Corollary 10.4.1]{Silverman2}
and results in the characterization given above.
\end{proof}

For any $\alpha \in \Q^{\times}$ and any $k \geq 2$ which is congruent to
$1$ modulo $h$, let
$$\psi^{k}_{\alpha} : \A_{K}^{\times}/K^{\times} \to K^{\times}$$
denote the Hecke character unramified away from $D$ and $\alpha$ such that
$$\psi^{k}_{\alpha}(\p) =
\psi(\p)^{(k-1)/h}\cdot \ls{\vep}{\p}_{w}^{(k-1)/h} \cdot
\ls{\alpha}{\p}_{w}$$
for any $\p$ relatively prime to $D$ and $\alpha$.
(The extra twist by $\vep$ will greatly
simplify the statements below.)  Note that $\psi^{k}_{\alpha}$ has
infinity type $(k-1,0)$.

\subsection{Modular forms}

Fix $k \geq 2$, $k \equiv 1 \amod{h}$, and $\alpha \in \Q^{\times}$.  
By \cite[Theorem 3.4]{Ribet} the Fourier series
$$g^{k}_{\alpha} := \sum_{(\a,D)=(\a,\alpha)=1} 
\psi^{k}_{\alpha}(\a)q^{\Nm \a}$$
(summing over ideals of $\O$ prime to $D$ and $\alpha$)
is a cusp form of weight $k$
which is an eigenform
for the Hecke operators $T_{n}$ with $n$ prime to $D$ and $\alpha$.  We let
$$f^{k}_{\alpha} = \sum a_{n}(f^{k}_{\alpha})q^{n}$$ 
denote the associated normalized newform; we have
$a_{n}(f_{\alpha}^{k}) = a_{n}(g_{\alpha}^{k})$ for $n$ prime to $D$
and $\alpha$.

We claim that $f^{k}_{\alpha}$ has rational Fourier coefficients.
Indeed, for $p$ relatively prime to $D$ and $\alpha$ 
Lemmas~\ref{lemma:gen} and~\ref{lemma:gen2} show that
$\psi^{k}_{\alpha}(\pbar) = \overline{\psi^{k}_{\alpha}(\p)}$,
so that
\begin{equation} \label{eq:conj}
a_{p}(f^{k}_{\alpha}) = \psi^{k}_{\alpha}(\p) +
\psi^{k}_{\alpha}(\pbar) \in \Q.
\end{equation}
By \cite[Corollary 3.1]{Ribet} the $a_{p}(f^{k}_{\alpha})$ for almost
all primes $p$ generate
the field of all Fourier coefficients of $f^{k}_{\alpha}$, so that
the rationality of $a_{n}(f^{k}_{\alpha})$ for all $n \geq 1$ follows.

Note that $f^{k}_{\alpha} = f^{k}_{-d\alpha}$,
but otherwise the $f^{k}_{\alpha}$ are all distinct.
In fact, any modular form of weight $k \equiv 1 \amod{h}$ with complex
multiplication by $K$ and with rational Fourier coefficients is equal to
$f^{k}_{\alpha}$ for some $\alpha \in \Q^{\times}$, 
although we will not prove this here.
(In general there may also be such forms of weights congruent to $1$ modulo
the exponent of the class group of $\O$, but these are the only other
possible weights.)

\begin{lemma} \label{lemma:ap}
Fix $m \geq 1$, $k \geq 2$, $k \equiv 1 \amod{h}$, and 
$\alpha \in \Q^{\times}$.
Let $p$ be a rational prime such that:
\begin{enumerate}
\item $p \equiv 1 \amod{m}$;
\item $p$ splits as $\p\pbar$ in $K/\Q$;
\item $p$ is relatively prime to $\alpha$.
\end{enumerate}
Then
$$\ls{a_{p}(f^{k}_{\alpha})}{p}_{m} = \ls{\psi^{k}_{\alpha}(\p)}{\pbar}_{m}.$$
\end{lemma}
\begin{proof}
Since $\pbar$ has norm $p$ we have by (\ref{eq:conj}) that
$$\ls{a_{p}(f^{k}_{\alpha})}{p}_{m} = 
\ls{\psi^{k}_{\alpha}(\p) + \psi^{k}_{\alpha}(\pbar)}{\pbar}_{m}.$$
Note that by Lemmas~\ref{lemma:gen} and~\ref{lemma:gen2}
$\psi(\pbar)$, and thus $\psi^{k}_{\alpha}(\pbar)$,
is divisible by $\pbar$.  Thus
$$\ls{\psi^{k}_{\alpha}(\p) + 
\psi^{k}_{\alpha}(\pbar)}{\pbar}_{m} = \ls{\psi^{k}_{\alpha}(\p)}{\pbar}_{m}$$
and the lemma follows.
\end{proof}

\subsection{Applications to abelian extensions of $K$}

In this section we assume for simplicity that $h=1$.
Fix a prime $m$ relatively prime to $w$ 
(that is, $m$ is odd and if $d=3$ then we also require $m \neq 3$)
and fix $p \equiv 1 \amod{m}$
which splits as $\p\pbar$ in $K/\Q$; assume also that $p > 3$.
Let $K_{p}^{m}$ denote the maximal abelian extension of $K$ of exponent
$m$ which is unramified away from $\p$ and $\pbar$.  

One can construct $K_{p}^{m}$ via the Hecke character
$\psi$ as follows.  Let
$\psi_{p} : G_{K} \to \O_{p}^{\times}$
be the $p$-adic Galois character associated to
$\psi$ via class field theory.  Let
$$\bar{\psi}_{p} : G_{K} \to (\O/p)^{\times} \cong \Fp^{\times} \times
\Fp^{\times}$$
be the reduction of $\tilde{\psi}_{p}$ 
(which is surjective by \cite[Corollay 5.20]{Rubin}) and let
$$\bar{\psi}_{p}^{m} : G_{K} \surj (\Z/m\Z)^{2}$$
denote the composition of $\bar{\psi}_{p}$ with
some fixed surjection $(\O/p)^{\times} \surj (\Z/m\Z)^{2}$.  We write
$K(\bar{\psi}_{p})$ and $K(\bar{\psi}_{p}^{m})$
for the fixed fields of the kernels of these characters.
Note that $K(\bar{\psi}_{p})$ is equal to $K(E[p])$ with $E$
a rational elliptic curve corresponding to $\psi$ over $K$.

\begin{lemma} \label{lemma:cft}
With notation as above,
$K_{p}^{m} = K(\bar{\psi}_{p}^{m})$.
\end{lemma}
\begin{proof}
By \cite[Theorem II.5.6]{Silverman2} the field
$K(\bar{\psi}_{p})$ contains the ray class field of $K$ of conductor $p$,
so that it certainly contains $K_{p}^{m}$.  Since
$K(\bar{\psi}_{p}^{m})$ 
is the maximal subextension of $K(\bar{\psi}_{p})/K$
of exponent $m$, to prove the lemma it thus suffices to show that
$K(\bar{\psi}_{p}^{m})/K$ is unramified away from $\p$ and $\pbar$.
To see this, note that $K(\bar{\psi}_{p})/K$ is unramified
away from $Dp$ since $\psi$ is unramified away from $D$.
Furthermore, by
\cite[Theorem 5.15]{Rubin} the image under $\psi_{p}$
of an inertia group at a prime dividing $D$ has order dividing $w$.
Since we are requiring $m$ to be relatively prime to $w$,
it follows that $K(\bar{\psi}_{p}^{m})/K$ must be 
unramified at all such primes, which completes the proof.
\end{proof}

It follows from Lemma~\ref{lemma:cft} that
$$\Gal(K_{p}^{m}/K) \cong (\Z/m\Z)^{2}.$$
In particular, $K_{p}^{m}$ is a Galois extension of $\Q$ of degree
$2m^{2}$; thus
$$e_{p}f_{p}g_{p} = 2m^{2}$$
where $e_{p}$ (resp.\ $f_{p}$, resp.\ $g_{p}$)
is the ramification index (resp.\ inertial degree, resp.\
splitting degree) of $p$ in $K_{p}^{m}/\Q$.
Since $K$ has no everywhere unramified extensions,
$K_{p}^{m}/\Q$ must be ramified at $p$;
since $p \neq m$, the inertia groups at $p$ must
be cyclic (as the tame inertia group of $\Qp$ is pro-cyclic), so that in fact
$e_{p}=m$.  As $2$ divides $g_{p}$, it follows that there are two
possibilities:
$$e_{p}=m,\qquad f_{p}=m, \qquad g_{p}=2$$
or
$$e_{p}=m,\qquad f_{p}=1, \qquad g_{p}=2m.$$
The next proposition shows that one can determine which of these occurs
in terms of the power residue of the $p^{\th}$ Fourier coefficient of $E$.

\begin{proposition} \label{prop:ravi}
Let $E$ be a rational elliptic curve with Hecke character $\psi$ over $K$.
For $p \equiv 1 \amod{m}$ greater than $3$ and split in $K/\Q$, we have
$f_{p}=1$ if and only if $\ls{a_{p}(E)}{p}_{m} = 1$.
\end{proposition}
\begin{proof}
Since
$$\Gal(K(E[p])/K) \cong (\O/p)^{\times} \cong (\Z/(p-1)\Z)^{2}$$
and
$$\Gal(K_{p}^{m}/K) \cong \Gal(K(E[p])/K)/m \cdot \Gal(K(E[p])/K),$$
one sees easily that $f_{p} = 1$ if and only if
the inertial degree of $\Qp(E[p])/\Qp$ divides $(p-1)/m$.
The latter condition is equivalent to the $p$-torsion of
$E$ over $\Fpbar$ being defined over ${\mathbf F}_{p^{(p-1)/m}}$.
Since $E$ is ordinary at $p$ (as $p$ splits in $K/\Q$) and thus has
one-dimensional $p$-torsion over $\Fpbar$, it follows that
to prove the proposition it suffices to show that
$\ls{a_{p}(E)}{p}_{m} = 1$ if and only if
$p$ divides $\# E({\mathbf F}_{p^{(p-1)/m}})$.

Let $\alpha,\beta$ be the $p$-adic roots of the Frobenius polynomial
$$x^{2} - a_{p}(E)x + p;$$
since $E$ is ordinary at $p$ we may order $\alpha,\beta$ so that
$\alpha$ is a $p$-adic unit and $\beta$ is
divisible by $p$.  In particular,
$$\alpha \equiv a_{p}(E) \pmod{p}.$$
By the Riemann hypotheses for elliptic curves over finite fields we have
\begin{align*}
\# E({\mathbf F}_{p^{(p-1)/m}}) &= p^{(p-1)/m} + 1 - \alpha^{(p-1)/m} -
\beta^{(p-1)/m} \\
&\equiv 1 - a_{p}(E)^{(p-1)/m} \pmod{p}.
\end{align*}
Since
$$a_{p}(E)^{(p-1)/m} \equiv \ls{a_{p}(E)}{p}_{m} \pmod{p},$$
the proposition follows from this.
\end{proof}

\begin{remark}
If $m$ is not relatively prime to $w$, then one can
easily obtain the analogue of Proposition~\ref{prop:ravi} via
Kummer theory.
\end{remark}

\begin{remark} \label{rmk:hwt}
Fix $k \geq 2$ and $\alpha \in \Q^{\times}$.  If $k-1$ is relatively
prime to $m$ (which in turn is still assumed relatively prime to $w$), then 
it follows from Lemma~\ref{lemma:ap} that
$$\ls{a_{p}(f^{k}_{\alpha})}{p}_{m} = \ls{a_{p}(E)}{p}_{m}.$$
Thus one can also compute $f_{p}$ via the Fourier coefficients of
$f^{k}_{\alpha}$.  
(In particular, this shows that it does not matter which rational
elliptic curve
with complex multiplication by $K$ one uses in Proposition~\ref{prop:ravi}.)
Unfortunately, the various hypotheses on $m$ and $k$ above
rule out any case in which we are able to calculate
$\ls{a_{p}(f^{k}_{\alpha})}{p}_{m}$, so that we not able to apply
Proposition~\ref{prop:ravi} in any generality.
\end{remark}

\begin{remark}
When considering Proposition~\ref{prop:ravi} for
$a_{p}(f^{k}_{\alpha})$ as above 
it is perhaps somewhat more
enlightening to regard $\Gal(K_{p}^{m}/K)$ as a quotient of the mod $p$ Galois
representation associated to $f^{k}_{\alpha}$; Proposition~\ref{prop:ravi}
can then be recovered from the fact that the restriction of this Galois
representation to a decomposition group at $p$ has the form
$$\left( \begin{array}{cc}
\chi^{k-1}\varphi^{-1} & * \\ 0 & \varphi \end{array} \right)$$
with $\chi$ the cyclotomic character and
$\varphi$ an unramified character with the
property that
$$\varphi(\Frob_{p}) \equiv a_{p}(f^{k}_{\alpha}) \pmod{p}.$$
\end{remark}

\section{Squares}

\subsection{Preliminaries}

In order to determine the quadratic character of Hecke
characters over imaginary quadratic fields with odd class number
we will need the following result.  Recall that
$\vep = -1$ (resp.\ $\vep = 1$) if $d \equiv 3 \amod{8}$
(resp.\ $d \equiv 7 \amod{8}$).

\begin{lemma} \label{lemma:pipibar}
Let $p$ be a prime which splits as $\p\pbar$ in $K/\Q$ and let $\pi$ be a
generator of $\p^{h}$.  (If $d=1$, further assume that $\pi \equiv
1 \amod{2+2i}$.)
\begin{enumerate}
\item If $p \equiv 1 \amod{4}$, then
$$\ls{\pi}{\pbar} = \ls{-d}{p}_{4}.$$
\item If $d \neq 2$ and $p \equiv 3 \amod{4}$, then
$$\ls{\pi}{\pbar} = \vep\ls{\pi}{\sqrt{-d}}.$$
\item If $d = 2$ and $p \equiv 3 \amod{4}$, then
$\ls{\pi}{\pbar} = 1$
if and only if $\pi$ is congruent to an element of (\ref{eq:2set})
modulo $4\sqrt{-2}$.
\end{enumerate}
\end{lemma}

\begin{proof}
Write $\pi = a + b\sd$ with $a,b \in
{\ts \oh \Z}$; if we write $a=2^{r}a'$, $b=2^{s}b'$ with $a',b'$ odd
integers, then $r=-1$ if and only if $s=-1$.  

We assume first that $p \equiv 1 \amod{4}$.  Since
$\pib \in \pbar$, we have
\begin{equation} \label{eq:p1m}
\ls{\pi}{\pbar} = \ls{\pi + \pib}{\pbar} = \ls{2b\sd}{\pbar} =
\ls{2b}{\pbar}\ls{\sd}{\pbar}.
\end{equation}
Since $2b$ is an integer and $\pbar$ has norm $p$, we have
$$\ls{2b}{\pbar} = \ls{2b}{p} = \ls{2}{p}^{s+1}\ls{b'}{p}.$$
By quadratic reciprocity and the fact that $h$ is odd we have
$$\ls{b'}{p} =
\ls{p}{b'} = \ls{p^{h}}{b'} = \ls{a^{2}+db^{2}}{b'} =
\ls{a^{2}}{b'} = 1.$$
As $\ls{\sd}{\p} = \ls{-d}{p}_{4}$
by (\ref{eq:p1m}) we conclude that
\begin{equation} \label{eq:seq}
\ls{\pi}{\pbar} = \ls{2}{p}^{s+1}\ls{-d}{p}_{4}.
\end{equation}
If $p \equiv 1 \amod{8}$ this completes the proof.  If $p \equiv 5 \amod{8}$,
then one shows easily (using that $\pi \equiv 1 \amod{2+2i}$ in the case $d=1$)
that $s= \pm 1$;
thus (\ref{eq:seq}) completes the proof in this case as well.

Next assume $p \equiv 3 \amod{4}$; in particular, we now have $d \neq 1$.
We have
$$\ls{\pi}{\pbar} = \ls{\pi+\pib}{\pbar}=\ls{2a}{\pbar} = \ls{2a}{p}
=\ls{2}{p}^{r+1}\ls{a'}{p}.$$
Using quadratic reciprocity and the fact that $h$ is odd, we obtain
\begin{multline} \label{eq:big}
\ls{\pi}{\pbar} = \ls{2}{p}^{r+1}\ls{-1}{a'}\ls{p}{a'}
= \ls{2}{p}^{r+1}\ls{-1}{a'}\ls{p^{h}}{a'} \\
= \ls{2}{p}^{r+1}\ls{-1}{a'}\ls{a^{2}+db^{2}}{a'}
= \ls{2}{p}^{r+1}\ls{-1}{a'}\ls{d}{a'}.
\end{multline}
For $d \neq 2$ we have $d \equiv 3 \amod{4}$, so that
a second application of quadratic reciprocity now yields
$$\ls{\pi}{\pbar} =
\ls{2}{p}^{r+1}\ls{a'}{d}
= \ls{2}{p}^{r+1}\ls{2}{d}^{r}\ls{a}{d}.$$
If $p \equiv d \amod{8}$, then this immediately yields
the lemma since $\ls{2}{p}=\ls{2}{d}=\vep$.
If $p \not \equiv d  \amod{8}$, one finds easily that
$r = \pm 1$; the lemma thus follows in this case as well.

When $d=2$, one must have $p \equiv 3 \amod{8}$ and $r=s=0$.
In particular, $\pi$ must be congruent to an element of
$$\{ 1+\sqrt{-2},3+\sqrt{-2},5+\sqrt{-2},7+\sqrt{-2},
1+3\sqrt{-2},3+3\sqrt{-2},5+3\sqrt{-2},7+3\sqrt{-2}\}$$
modulo $4\sqrt{-2}$.  By (\ref{eq:big}) we have
$$\ls{\pi}{\pbar} = \ls{2}{p}\ls{-2}{a} = -\ls{-2}{a},$$
which is $1$ if and only if $a=5,7$.  The lemma thus follows
from the definition of (\ref{eq:2set}).
\end{proof}

\subsection{Densities}

Combining Lemma~\ref{lemma:ap} with
Lemma~\ref{lemma:pipibar} 
yields the following result on the quadratic character
of the modular forms $f^{k}_{\alpha}$.

\begin{proposition} \label{prop:sqr}
Fix $k \geq 2$, $k \equiv 1 \amod{h}$, and $\alpha \in \Q^{\times}$.
Let $p$ be a rational prime,
relatively prime to $\alpha$, which splits as $\p\pbar$ in $K/\Q$.
\begin{enumerate}
\item For $d \neq 1$,
$$\ls{a_{p}(f^{k}_{\alpha})}{p} =
\begin{cases} 1 & p \equiv 1 \amod{4} \text{~and $k$ odd}; \\
 \ls{-d}{p}_{4} & p \equiv 1 \amod{4} \text{~and $k$ even}; \\
\ls{\alpha}{p} & p \equiv 3 \amod{4}.
\end{cases}$$
\item For $d = 1$,
$$\ls{a_{p}(f^{k}_{\alpha})}{p} =
\begin{cases} 1 & p \equiv 1 \amod{8}; \\
\ls{\alpha}{p} & p \equiv 5 \amod{8}.
\end{cases}$$
\end{enumerate}
\end{proposition}
\begin{proof}
Assume first that $d \neq 1,3$, so that $w=2$.
By Lemma~\ref{lemma:ap} and the definition of $\psi^{k}_{\alpha}$ we have
\begin{multline*}
\ls{a_{p}(f^{k}_{\alpha})}{p} =
\ls{\psi(\p)^{(k-1)/h} \cdot \ls{\vep}{\p}^{(k-1)/h} 
\cdot \ls{\alpha}{\p}}{\pbar} \\
= \ls{\psi(\p)}{\pbar}^{(k-1)/h} \cdot \ls{\ls{\vep}{p}}{p}^{(k-1)/h} \cdot
\ls{\ls{\alpha}{p}}{p}.
\end{multline*}
When $p \equiv 1 \amod{4}$, this yields
$$\ls{a_{p}(f^{k}_{\alpha})}{p} = \ls{\psi(\p)}{\pbar}^{(k-1)/h} =
\ls{-d}{p}_{4}^{(k-1)/h}$$
by Lemma~\ref{lemma:pipibar} and the fact that $\psi(\p)$ generates
$\p^{h}$.  Since $\ls{-d}{p}_{4} = \pm 1$ (as $p$ splits in
$K$) and $(k-1)/h$ is even if and only if $k$ is odd, the lemma follows
in this case.

When $p \equiv 3 \amod{4}$
we instead obtain
$$\ls{a_{p}(f^{k}_{\alpha})}{p} = \ls{\psi(\p)}{\pbar}^{(k-1)/h} \cdot
\vep^{(k-1)/h} \cdot \ls{\alpha}{p}.$$
By Lemmas~\ref{lemma:gen}, \ref{lemma:gen2} and~\ref{lemma:pipibar} we have
$\ls{\psi(\p)}{\pbar} = \vep$, so that the lemma follows in this case as
well.

The proof for $d=3$ is similar using that
$$\ls{\ls{\alpha}{p}_{6}}{\pbar} = \begin{cases} 1 & p \equiv 1 \amod{12}; \\
\ls{\alpha}{p} & p \equiv 7 \amod{12}. \end{cases}$$
The proof for $d=1$ also proceeds similarly, taking into account that
$\vep = 2$, that $\ls{-1}{p}_{4}=1$ for
$p \equiv 1 \amod{8}$, and that
$$\ls{\ls{\alpha}{p}_{4}}{\pbar} = \begin{cases} 1 & p \equiv 1 \amod{8}; \\
\ls{\alpha}{p} & p \equiv 5 \amod{8}. \end{cases}$$
\end{proof}

It is now a simple matter to determine the density of squares among the
non-zero Fourier coefficients of $f^{k}_{\alpha}$.
By the definition of $f^{k}_{\alpha}$ 
we see that $a_{p}(f^{k}_{\alpha}) = 0$ for $p$ inert
in $K/\Q$ while $p \nmid a_{p}(f^{k}_{\alpha})$ 
for $p$ split in $K/\Q$ and relatively prime to $\alpha$.
Thus
\begin{align*}
\delta_{2}(f^{k}_{\alpha}) &:=
\frac{\text{density of~} \left\{p \,;\, 
\ls{a_{p}(f^{k}_{\alpha})}{p} = 1 \right\}}{\text{density of~} \{p \,;\,
a_{p}(f^{k}_{\alpha}) \not\equiv 0 \amod{p} \}} \\ &=
\frac{\text{density of~} \left\{p \,;\, 
\ls{a_{p}(f^{k}_{\alpha})}{p} = 1 \right\}}{\text{density of~} \{p \,;\,
p \text{~split in~} K/\Q \}} \\ &= 2 \cdot \left(
\text{density of~} \left\{p \,;\, 
\ls{a_{p}(f^{k}_{\alpha})}{p} = 1 \right\} \right)
\end{align*}
if this density is defined.

\begin{theorem} \label{thm:sqrs}
Fix $k \geq 2$, $k \equiv 1 \amod{h}$, and $\alpha \in \Q^{\times}$.
\begin{itemize}
\item For $d \neq 1$ and $k$ even
$$\delta_{2}(f^{k}_{\alpha}) = \begin{cases}
\frac{3}{4} & D(\alpha) \in \{1,D\}; \\
\frac{1}{4} & D(\alpha) \in \{-4,4d\}; \\
\frac{1}{2} & \text{otherwise.}
\end{cases}$$
\item For $d \neq 1$ and $k$ odd
$$\delta_{2}(f^{k}_{\alpha}) = \begin{cases}
1 & D(\alpha) \in \{1,D\}; \\
\frac{1}{2} & D(\alpha) \in \{-4,4d\}; \\
\frac{3}{4} & \text{otherwise}.
\end{cases}$$
\item For $d = 1$,
$$\delta_{2}(f^{k}_{\alpha}) = \begin{cases}
1 & D(\alpha) \in \{1,-4\}; \\
\frac{1}{2} & D(\alpha) \in \{\pm 8\}; \\
\frac{3}{4} & \text{otherwise}.
\end{cases}$$
\end{itemize}
\end{theorem}
\begin{proof}
We consider $d \neq 1$; the proof for $d=1$ is handled in an entirely similar
fashion, taking into account the different form
of Proposition~\ref{prop:sqr} in this case.
Consider first primes $p \equiv 1 \amod{4}$, $p$ relatively prime to $\alpha$ 
and split in $K/\Q$.  These primes are (up to a finite set) precisely
those which split completely in $\Q(\sqrt{-d},i)/\Q$.
By Proposition~\ref{prop:sqr}, if $k$ is odd then every such $p$
satisfies $\ls{a_{p}(f^{k}_{\alpha})}{p} = 1$, while if $k$
is even such a $p$ satisfies $\ls{a_{p}(f^{k}_{\alpha})}{p} = 1$ if and only if
it splits completely in
$\Q(\sqrt[4]{-d},i)/\Q$.
By the Cebatorev density theorem the set of such primes
has density $\frac{1}{4}$ or $\frac{1}{8}$, respectively.

Next consider primes $p \equiv 3 \amod{4}$, $p$ relatively prime to $\alpha$
and split in $K/\Q$.  These primes are inert in $\Q(i)/\Q$, and by
Proposition~\ref{prop:sqr} they satisfy $\ls{a_{p}(f^{k}_{\alpha})}
{p}=1$ if and
only if they split in $\Q(\sqrt{\alpha})/\Q$.  In particular, the set
of such primes has density $\frac{1}{8}$ unless
$\Q(\sqrt{\alpha})$ lies in $\Q(\sqrt{-d},i)$.
This occurs only if
$D(\alpha) \in \{1,D\}$ or $D(\alpha') \in \{-4,4d\}$, in which case
these primes have density $\frac{1}{4}$ or $0$, respectively.
The proposition now follows easily in these cases.
\end{proof}

\section{Higher powers}

\subsection{Cubes}

For the field $K=\Q(\sqrt{-3})$
it is possible to use the law of cubic reciprocity to study the
cubic character of the Fourier coefficients of the modular form
$f^{k}_{\alpha}$.

\begin{proposition} \label{prop:cube}
Fix $k \geq 2$ and $\alpha \in \Q^{\times}$.
Let $p \equiv 1 \amod{3}$ be a rational prime,
relatively prime to $\alpha$, which splits in $K/\Q$.  Then
$$\ls{a_{p}(f^{k}_{\alpha})}{p}_{3} =
\begin{cases} 1 & p \equiv 1 \amod{9}; \\
\ls{\alpha}{p}_{3}^{2} & p \equiv 4 \amod{9}; \\
\ls{\alpha}{p}_{3} & p \equiv 7 \amod{9}.
\end{cases}$$
\end{proposition}
\begin{proof}
Let $\pi$ be a prime divisor of $p$ which is congruent to $1$ modulo
$3$.
By the law of cubic reciprocity (see, for example,
\cite[Section 7.2]{Lemmermeyer}) we have
$$\ls{\pi}{\pib}_{3} = \ls{\pib}{\pi}_{3}.$$
On the other hand, 
$$\overline{\ls{\pi}{\pib}}_{3} = \ls{\pib}{\pi}_{3};$$
it follows that in fact $\ls{\pi}{\pib}_{3} = 1$.
In particular, by Lemma~\ref{lemma:ap} and the
definition of $\psi^{k}_{\alpha}$ we have
\begin{multline*}
\ls{a_{p}(f^{k}_{\alpha})}{p}_{3} =
\ls{\pi}{\pib}_{3}^{k-1}\ls{\ls{-1}{p}}{p}_{3}^{k-1}
\ls{\ls{\alpha}{p}_{6}}{\pib}_{3} =
\ls{\ls{\alpha}{p}_{6}}{\pib}_{3} \\
= \ls{\alpha}{p}_{6}^{(p-1)/3} =
\ls{\alpha}{p}_{3}^{(p-1)/6};
\end{multline*}
the proposition follows easily.
\end{proof}

\begin{theorem} \label{thm:cube}
Fix $k \geq 2$ and $\alpha \in \Q^{\times}$.  Then
$$\delta_{3}(f^{k}_{\alpha}) = \begin{cases}
1 & \alpha \in \Q^{\times 3}; \\
\frac{5}{9} & \alpha \notin \Q^{\times 3}.
\end{cases}$$
\end{theorem}
\begin{proof}
If $\alpha \in \Q^{\times 3}$, then $\ls{\alpha}{p}_{3} = 1$ for
all $p \equiv 1 \amod{3}$, so that $\delta_{3}(f_{\psi}) = 1$
by Proposition~\ref{prop:cube}.  On the other hand, if
$\alpha \notin \Q^{\times 3}$, then primes $p \equiv 1 \amod{9}$ always
satisfy $\ls{a_{p}(\psi)}{p}_{3} = 1$, while primes $p \equiv 4,7 \amod{9}$
satisfy $\ls{a_{p}(\psi)}{p}_{3} = 1$ if and only if
$\ls{\alpha}{p}_{3} = 1$.  The former condition is
equivalent to $p$ splitting completely in $\Q(\mu_{9})/\Q$, so that
the set of such $p$ has density $\frac{1}{6}$.
The latter condition is equivalent to $p$ having inertial degree
$3$ in $\Q(\mu_{9})/\Q$ and splitting completely in
$\Q(\mu_{9},\sqrt[3]{\alpha})/\Q(\mu_{9})$.  
The unique cubic subfield
of $\Q(\mu_{9})$ is $\Q(\mu_{9})^{+}$, which is Galois over $\Q$ and
thus can not contain any non-trivial cube roots.  Thus
$\Q(\mu_{9}) \cap \Q(\sqrt[3]{\alpha}) = \Q$, so that
the set of such $p$ has
density $\frac{2}{3} \cdot \frac{1}{2} 
\cdot \frac{1}{3} = \frac{1}{9}$ by the Cebatorev
theorem.  
$$
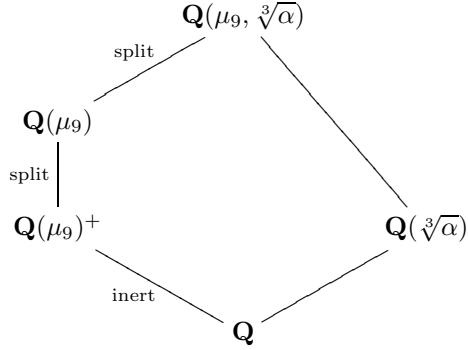
\xymatrix{
& {\Q(\mu_{9},\sqrt[3]{\alpha})} \ar@{-}[dl]_{\text{split}} \ar@{-}[ddr] & \\
{\Q(\mu_{9})} \ar@{-}[d]_{\text{split}} & & \\
{\Q(\mu_{9})^{+}} \ar@{-}[dr]_{\text{inert}} 
& & {\Q(\sqrt[3]{\alpha})} \ar@{-}[dl] \\
& {\Q} &}$$
Combining these results, we see that the set of $p \equiv 1 \amod{3}$
with $\ls{a_{p}(f^{k}_{\alpha})}{p}_{3} = 1$ 
has density $\frac{5}{18}$, from which
it follows that $\delta_{3}(f^{k}_{\alpha}) = \frac{5}{9}$, as claimed.
\end{proof}

\subsection{$m^{\th}$ powers, $m | k-1$}

For any $\alpha \in \Q^{\times}$ we write $D(\alpha)$ for the discriminant of
the quadratic field $\Q(\sqrt{\alpha})$ over $\Q$.  For simplicity
we state the next result only for $d > 3$; the cases become
overwhelming for $d \leq 3$.

\begin{theorem} \label{prop:high}
Assume that $d > 3$ and fix $k \geq 2$, $k \equiv 1 \amod{h}$, 
and $\alpha \in \Q^{\times}$.  Then for any $m|k-1$ we have
$$\delta_{m}(f^{k}_{\alpha}) = \begin{cases}
1 & m \text{~odd or~} D(\alpha) | m \text{~or~}
D(-d\alpha) | m; \\
\frac{1}{2} & m \text{~even and~} D(\alpha) \nmid m \text{~and~}
D(-d\alpha) \nmid m \text{~and either~}\\ 
& D(\alpha) | 2m \text{~or~} D(-d\alpha) | 2m; \\
\frac{3}{4} & otherwise.
\end{cases}$$
\end{theorem}
\begin{proof}
Let $p \equiv 1 \amod{m}$ be a prime which splits as $\p\pbar$ in $K/\Q$.
By Lemma~\ref{lemma:ap} and the definition of $\psi_{\alpha}^{k}$ we have
$$\ls{a_{p}(f^{k}_{\alpha})}{p}_{m} = 
\ls{\psi'{}^{k-1}(\PP) \cdot \ls{\vep}{p}^{(k-1)/h}
\cdot \ls{\alpha}{p}}{\pbar}_{m} = 
\ls{\ls{\vep}{p}^{(k-1)/h} \cdot \ls{\alpha}{p}}{p}_{m}$$
for any prime $\PP$ of $H$ lying over $\p$.
Since $h$ is odd and $\ls{\vep}{p} = \pm 1$, we certainly have
$$\ls{\ls{\vep}{p}}{p}_{m}^{(k-1)/h} = 1,$$
so that
$$\ls{a_{p}(f^{k}_{\alpha})}{p}_{m} = 
\ls{\ls{\alpha}{p}}{p}_{m}.$$
This is not equal to $1$ if and only if
$\ls{\alpha}{p} = \ls{-1}{p}_{m} = -1$.  (Note that $\ls{-1}{p}_{m}
= \pm 1$ since $p \equiv 1 \amod{m}$.)
This in turn is equivalent to the following splitting behavior of $p$:
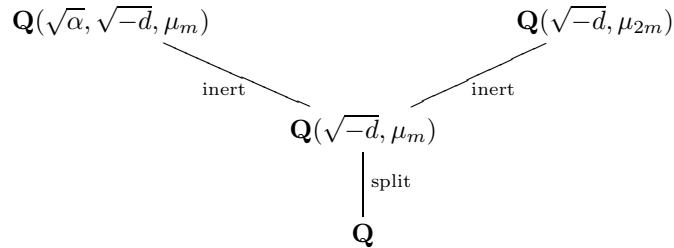
$$\xymatrix{
{\Q(\sqrt{\alpha},\sqrt{-d},\mu_{m})} \ar@{-}[dr]_{\text{inert}} & &
{\Q(\sqrt{-d},\mu_{2m})} \ar@{-}[dl]^{\text{inert}} \\
& {\Q(\sqrt{-d},\mu_{m})} \ar@{-}[d]^{\text{split}} & \\
& {\Q} &}$$

In particular, there are no such primes (so that $\delta_{m}
(f^{k}_{\alpha})=1$)
if and only if either of the top two extensions are trivial.
If both of these extensions
are non-trivial, then it follows from the Cebatorev theorem that
$\delta_{m}(f^{k}_{\alpha}) = \frac{3}{4}$ unless the two extensions
coincide, in which case $\delta_{m}(f^{k}_{\alpha}) = \frac{1}{2}$.
One now checks easily using Lemma~\ref{lemma:cyc} below
(and the fact that $\Q(\mu_{m}) = \Q(\mu_{2m})$ if and only if $m$ is odd)
that the conditions given in the statement of the theorem are equivalent 
to these field theoretic conditions.
\end{proof}

We remark that this result recovers
Theorem~\ref{thm:sqrs} when $k$ is odd.

\begin{lemma} \label{lemma:cyc}
Fix $n \geq 1$ and $\alpha,\beta \in \Q^{\times}$.
Then
$$[\Q(\sqrt{\alpha},\sqrt{\beta},\mu_{n}):\Q] =
\begin{cases}
\varphi(n) & D(\alpha) \text{~and~} D(\beta) \text{~divide~}n; \\
2\varphi(n) & \text{exactly one of~} D(\alpha),D(\beta),D(\alpha\beta)
\text{~divides~}n; \\
4\varphi(n) & \text{otherwise.}
\end{cases}$$
\end{lemma}
\begin{proof}
Set $F=\Q(\sqrt{\alpha},\sqrt{\beta})$.
We have
$$[F(\mu_{n}):\Q] =
[F:F \cap \Q(\mu_{n})] \cdot [\Q(\mu_{n}):\Q] =
[F:F \cap \Q(\mu_{n})] \cdot \varphi(n).$$
Note that $\Q(\mu_{D(\alpha)})$ is the smallest cyclotomic field
containing $\sqrt{\alpha}$, so that in general $\sqrt{\alpha} \in \Q(\mu_{n})$
if and only if $D(\alpha)$ divides $n$.  As every subfield of
$F$ is generated by some subset of $\{ \sqrt{\alpha},\sqrt{\beta},\sqrt{\alpha
\beta}\}$, this allows one to easily compute $F \cap \Q(\mu_{n})$
in terms of $D(\alpha)$, $D(\beta)$, and $D(\alpha\beta)$.
The lemma follows from this computation.
\end{proof}


\begin{thebibliography}{10}

\bibitem{Cremona}
John Cremona, \emph{Algorithms for modular elliptic curves, Second Edition}.
Cambridge University Press, Cambridge, 1997.

\bibitem{Deligne}
Pierre Deligne, \emph{La conjecture de Weil. I.},
Inst.\ Hautes \'Etudes Sci. {\bf 
43} (1974), 273--307.

\bibitem{deShalit}
Ehud de Shalit, \emph{Iwasawa theory of elliptic curves with complex 
multiplication}, Academic Press, Boston, MA, 1987.

\bibitem{FM}
Jean-Marc Fontaine and Barry Mazur,
\emph{Geometric Galois representations}, in:
\emph{Elliptic curves, modular forms and Fermat's last theorem
(Hong Kong, 1993)}, International Press, Cambridge, MA, 1995, 
pp.\ 41--78.

\bibitem{Gross}
Benedict Gross, \emph{Arithmetic on elliptic curves with complex
multiplication}.  Lecture Notes in Math.\ {\bf 776},
Springer--Verlag, Berlin, 1980.

\bibitem{Ravi}
Chandrashekhar Khare, Michael Larsen and Ravi Ramakrishna,
\emph{Construction of semisimple $p$-adic Galois representations with
prescribed properties}, preprint.

\bibitem{Lemmermeyer}
Franz Lemmermeyer, \emph{Reciprocity laws: From Euler to Eisenstein}.
Springer--Verlag, Berlin, 2000.

\bibitem{Ribet}
Ken Ribet, \emph{Galois representations attached to eigenforms with
nebentypus}, in: \emph{Modular functions of one variable, V},
Lectures Notes in Math.\ {\bf 601}, Springer--Verlag, Berlin, 1977,
pp.\ 7--51.

\bibitem{Rubin}
Karl Rubin, \emph{Elliptic curves with complex multiplication and the
conjecture of Birch and Swinnerton-Dyer}, in:
\emph{Arithmetic theory of elliptic curves}, Lectures Notes in
Math.\ {\bf 1716}, Springer--Verlag, Berlin, 1999, pp.\ 167--234.

\bibitem{Scholl}
Anthony Scholl, \emph{Motives for modular forms},
Inventiones Math.\ {\bf 100} (1990), 419--430.

\bibitem{Silverman2}
Joseph Silverman, \emph{Advanced topics in the arithmetic of elliptic curves}.
Graduate Texts in Math.\ {\bf 151}, Springer--Verlag, New York, 1994.

\bibitem{Stein}
William Stein, \emph{The modular forms explorer},
available at:\\ \texttt{http://modular.fas.harvard.edu/mfd/mfe/}.

\end{thebibliography}
\end{document}